\documentclass[a4paper,11pt,twoside,reqno]{amsart}

\usepackage[utf8]{inputenc}
\usepackage[plainpages=false,pdfpagelabels=true]{hyperref}
\usepackage{amssymb,amsthm}
\usepackage[margin=1in]{geometry}
\usepackage[all]{xy}
\newtheorem{Satz}{Theorem}[section]
\newtheorem{Prop}[Satz]{Proposition}
\newtheorem{Lem}[Satz]{Lemma}

\newtheorem{Cor}[Satz]{Corollary}
\newcommand{\vol}{{\operatorname{Vol}}}
\theoremstyle{definition}
\newtheorem{Dfn}[Satz]{Definition}
\newtheorem{Bem}[Satz]{Remark}
\newtheorem{Bsp}[Satz]{Example}
\allowdisplaybreaks[1]

\renewcommand{\epsilon}{\varepsilon}

\newcommand{\R}{\ensuremath{\mathbb{R}}}

\newcommand{\Z}{\ensuremath{\mathbb{Z}}}
\newcommand{\N}{\ensuremath{\mathbb{N}}}

\numberwithin{equation}{section}

\title{An estimate on the nodal set of eigenspinors on closed surfaces}
\author{Volker Branding}
\date{\today}
\address{University of Vienna, Faculty of Mathematics\\
Oskar-Morgenstern-Platz 1, 1090 Vienna, Austria}
\email[]{volker.branding@univie.ac.at}
\begin{document}

\begin{abstract}
We use a modified Bochner technique to derive an inequality relating 
the nodal set of eigenspinors to eigenvalues of the Dirac operator on closed surfaces.
In addition, we apply this technique to solutions of similar spinorial equations.
\end{abstract} 

\maketitle

\section{Introduction and results}
Throughout this note we assume that \((M,g)\) is a closed, connected, oriented surface
with a Riemannian metric \(g\). The bundle \(SO(M)\) of oriented orthonormal
frames is an \(S^1\)-principal bundle over the surface \(M\).
Let \(\Theta:S^1\to S^1\) be the nontrivial double covering of \(S^1\).
A spin structure on \(M\) is an \(S^1\)-principal bundle \(Spin(M)\) over \(M\) 
together with a twofold covering map \(\theta\colon Spin(M)\to SO(M)\) such that the diagram
\begin{center}
$$
\xymatrix{ 
\text{Spin}(M)\times S^1 \ar[dd]^{\theta\times\Theta} \ar[r] & \text{Spin}(M) \ar[dd]^{\theta} \ar[dr] &\\
& & M \\
\text{SO}(M)\times S^1\ar[r] & \text{SO}(M) \ar[ur]&}
$$
\end{center}
commutes. Every orientable surface admits a spin structure, 
the number of possible spin structure is equal to the number of elements in \(H^1(M,\Z_2)\).
On the spinor bundle \(\Sigma M\) we have a metric connection \(\nabla\) and a hermitian scalar product 
of which we will always take the real part turning it into an euclidean scalar product.
Sections in the spinor bundle are called \emph{spinors}.
Moreover, there exists the canonical splitting of the spinor bundle \(\Sigma M\) into the \emph{bundle of positive spinors} \(\Sigma^+M\)
and the \emph{bundle of negative spinors} \(\Sigma^-M\), that is \(\Sigma M=\Sigma^+M\oplus\Sigma^-M\).

In addition, we have the Clifford multiplication of spinors with tangent vectors, denoted by \(X\cdot\psi\) for
\(X\in TM\) and \(\psi\in\Gamma(\Sigma M)\). Clifford multiplication is skew-symmetric
\[
\langle X\cdot\psi,\xi\rangle=-\langle \psi,X\cdot\xi\rangle
\]
and satisfies the Clifford relations
\[
X\cdot Y\cdot\psi+Y\cdot X\cdot\psi=-2g(X,Y)\psi
\]
for \(X,Y\in TM\) and \(\psi,\xi\in\Gamma(\Sigma M)\).

\begin{Dfn}
The \emph{Dirac operator} \(D\) maps smooth sections of \(\Sigma M\) to smooth
sections of \(\Sigma M\) and is given by
\[
D\psi:=e_1\cdot\nabla_{e_1}\psi+e_2\cdot\nabla_{e_2}\psi,
\]
where \(e_1,e_2\) is a local orthonormal frame of \(TM\).
\end{Dfn}
The Dirac operator is a first order, elliptic operator, which is self-adjoint with respect to the \(L^2\)-norm.
Thus, from general spectral theory we know that the spectrum of the Dirac operator is real, discrete
and tends rapidly to infinity. In addition, it is well known that in contrast to the Laplacian
the spectrum of the Dirac operator consists of both positive and negative eigenvalues.

The square of the Dirac operator satisfies the \emph{Schrödinger-Lichnerowicz} formula
\begin{equation}
\label{schroedinger-lichnerowicz}
D^2=\nabla^*\nabla+\frac{R}{4},
\end{equation}
where \(R\) denotes the scalar curvature of the manifold.

For more background material on spin geometry and the Dirac operator we refer the reader to the books \cite{MR1031992} and \cite{MR1476425}.

\begin{Dfn}
A spinor \(\psi\in\Gamma(\Sigma M)\) is called eigenspinor with eigenvalue \(\lambda\) if it satisfies
\begin{equation}
D\psi=\lambda\psi.
\end{equation}
In particular, an eigenspinor corresponding to the eigenvalue \(\lambda=0\) is called \emph{harmonic}.
\end{Dfn}
In general, the spectrum of the Dirac operator \(D\) cannot be computed explicitly.
There are only few manifolds with high symmetry, who allow to explicitly determine the spectrum,
for example flat tori  \cite{MR750754} and round spheres \cite{MR1361548}.

However, it is possible to estimate the spectrum. 
A fundamental inequality for the eigenvalues of the Dirac operator is Friedrich's inequality \cite{MR600828}
\begin{equation}
\label{friedrich-inequal}
\lambda^2\geq\frac{n}{4(n-1)}\inf_M R,
\end{equation}
where \(n\) is the dimension of the manifold  \(M\).

On closed surfaces, the following inequality was given by Bär in \cite{MR1162671}
\begin{equation}
\label{baer-inequal}
\lambda^2\geq\frac{2\pi\chi(M)}{\vol(M,g)},
\end{equation}
where \(\chi(M)\) is the Euler characteristic of \(M\).
It directly follows from \eqref{baer-inequal} that there do not exist harmonic spinors on \(S^2\). 
However, on surfaces of genus \(g\geq 1\) there always exist a metric and a spin structure 
admitting harmonic spinors, see \cite{MR0358873}, Proposition 2.4, \cite{MR2509837}, Theorem 6.2.1 and \cite{MR1186015}.
It is also possible to give estimates on the first non-zero eigenvalue of the Dirac operator that include 
the spin structure, see for example \cite{MR1900319}.

For more details on the spectrum of the Dirac operator see the book \cite{MR2509837}.

Basically, inequality \eqref{baer-inequal} gives information about small eigenvalues of the Dirac operator.
In this note we will derive an inequality that provides properties of large eigenvalues.
On closed surfaces the nodal set of eigenspinors is discrete \cite{MR1714341},
which enables us to prove the following 
\begin{Satz}
\label{main-result}
Suppose that \((M,g)\) is a closed, connected spin surface with fixed spin structure. Then every eigenvalue \(\lambda\)
of the Dirac operator \(D\) satisfies the following inequality
\begin{equation}
\label{eigenvalue-equation-improved}
\lambda^2\geq\frac{2\pi\chi(M)}{\vol(M,g)}+\frac{4\pi N(\psi)}{\vol(M,g)},
\end{equation}
where \(\chi(M)\) is the Euler characteristic of \(M\). Moreover, \(N(\psi)\) denotes the sum of the order of the zero's of the corresponding eigenspinor \(\psi\),
that is
\begin{equation}
N(\psi)=\sum_{p\in M,|\psi|(p)=0}n_p.
\end{equation}
\end{Satz}

\begin{Cor}
Of course, \eqref{eigenvalue-equation-improved} can also be interpreted as an inequality for the nodal set of an eigenspinor \(\psi\)
of the Dirac operator
\begin{equation}
\label{estimate-nodal-set}
N(\psi)\leq\frac{\vol(M,g)\lambda^2}{4\pi}-\frac{\chi(M)}{2}.
\end{equation}
\end{Cor}

We will discuss a similar inequality for harmonic spinors, twistor spinors and solutions of
a semi-linear Dirac equation.

\begin{Bem}
Note that \eqref{estimate-nodal-set} also holds if \(\lambda\) is an eigenvalue with higher multiplicity.
In this case we could consider a linear combination of eigenspinors and the nodal set could depend on the 
particular linear combination. However, the estimate \eqref{estimate-nodal-set} holds for all linear combinations.
\end{Bem}

\begin{Bem}
We cannot expect that a similar inequality holds in higher dimensions since
the nodal set of eigenspinors is no longer discrete \cite{MR1714341}.
\end{Bem}

\section{Proof of the main Theorem}
By the main result of \cite{MR1714341} we know that on a two-dimensional manifold the zero-set of
eigenspinors is discrete. In the following we will make use of the energy-momentum tensor \(T(X,Y)\)
\begin{equation}
T(X,Y):=\langle X\cdot\nabla_Y\psi+Y\cdot\nabla_X\psi,\psi\rangle.
\end{equation}
This tensor arises if one varies the functional \(E(\psi)=\int_M\langle\psi,D\psi\rangle dM\) with respect to the metric.

The following Lemma can be found in \cite{MR1738150}, Lemma 5.1, see also \cite{MR1806438}.
Since we need a slightly more general version of it we also give a proof here.
In contrast to the reference given above we use the analyst's sign convention for the Laplacian.
\begin{Lem}
\label{lem-modified-connection}
Suppose that \(\psi\in\Gamma(\Sigma M)\) does not have any zeros. Then the following inequality holds
\begin{equation}
\label{eigenvalue-modified-connection}
\frac{\langle\psi,D^2\psi\rangle}{|\psi|^2}\geq\frac{R}{4}+\frac{|T|^2}{4|\psi|^4}-\Delta\log|\psi|
-\langle D\psi,d(\log|\psi|^2)\cdot\psi\rangle.
\end{equation}
\end{Lem}
\begin{proof}
We set
\[
\tilde{\nabla}_X\psi:=\nabla_X\psi-2\alpha(X)\psi-\beta(X)\cdot\psi-X\cdot\alpha\cdot\psi, 
\]
with a one-form \(\alpha\) and a symmetric \((1,1)\)-tensor \(\beta\) given by
\[
\alpha:=\frac{d|\psi|^2}{2|\psi|^2}, \qquad \beta:=-\frac{T(\cdot,\cdot)}{2|\psi|^2}.
\]
By a direct computation using the real scalar product on \(\Sigma M\) we find summing over repeated indices
\begin{align*}
|\tilde{\nabla}\psi|^2=|\nabla\psi|^2+2|\alpha|^2|\psi|^2+|\beta|^2|\psi|^2
-4\alpha({e_i})\langle\nabla_{e_i}\psi,\psi\rangle+2\langle\beta(e_i)\cdot\nabla_{e_i}\psi,\psi\rangle
+2\langle D\psi,\alpha\cdot\psi\rangle.
\end{align*}
Moreover, we have
\begin{align*}
|\nabla\psi|^2=&\langle\psi,D^2\psi\rangle-\frac{R}{4}|\psi|^2+\frac{1}{2}\Delta|\psi|^2, \\
\alpha({e_i})\langle\nabla_{e_i}\psi,\psi\rangle=&|\alpha|^2|\psi|^2=\frac{|d|\psi|^2|^2}{4|\psi|^2},\\
\langle\beta(e_i)\cdot\nabla_{e_i}\psi,\psi\rangle=&-\frac{|T|^2}{4|\psi|^2}.
\end{align*}
Thus, we arrive at
\[
0\leq|\tilde{\nabla}\psi|^2=\langle\psi,D^2\psi\rangle-\frac{R}{4}|\psi|^2+\frac{1}{2}\Delta|\psi|^2-\frac{|d|\psi|^2|^2}{2|\psi|^2}
-\frac{|T|^2}{4|\psi|^2}+2\langle D\psi,\alpha\cdot\psi\rangle
\]
yielding the result.
\end{proof}

The following Lemma will be the key-tool for the further analysis.
For the sake of completeness we also present a proof, where we follow \cite{MR1474501}.

\begin{Lem}
Suppose \(M\) is a closed Riemannian surface. If the zero set of \(|\psi|\) is discrete and \(|\psi|\) does not vanish identically, 
then the following equality holds
\begin{equation}
\label{laplace-log}
\int_M\Delta\log|\psi|dM=-2\pi\sum_{p\in M,|\psi|(p)=0}n_p,
\end{equation}
where \(n_p\) is the order of \(|\psi|\) at the point \(p\).
\end{Lem}

\begin{proof}
Since the zeros of \(|\psi|\) are isolated and \(M\) is compact,
the number of zeros \(p_1,\ldots,p_k\) is finite in \(M\).
Let \(D_\epsilon(p_j)\) be a small disc of radius \(\epsilon>0\) around \(p_j\).
Applying the divergence theorem we get
\begin{align*}
\int_M\Delta\log|\psi|dM=\lim_{\epsilon\to 0}\int_{M\setminus\cup_{j=1}^k D_\epsilon(p_j)}\Delta\log|\psi|dM
=-\lim_{\epsilon\to 0}\sum_{j=1}^k\int_{\partial D_\epsilon(p_j)}\frac{\partial}{\partial r}\log|\psi|d\theta, 
\end{align*}
where \(\frac{\partial}{\partial r}\) denotes the radial derivative.
Using a local Taylor expansion for each \(j\) we get
\[
\int_{\partial D_\epsilon(p_j)}\frac{\partial}{\partial r}\log|\psi|d\theta=2\pi n_{p_j}+O(\epsilon),
\]
where \(n_j\) denotes the order of the first non-vanishing term in the Taylor expansion of \(\psi\) in \(p_j\).
By letting \(\epsilon\to 0\) we thus obtain
\[
\int_M\Delta\log|\psi|dM=-2\pi\sum_{j=1}^kn_{p_j}=-2\pi\sum_{p\in M,|\psi|(p)=0}n_p,
\]
which proves the Lemma.
\end{proof}

Finally, we apply Lemma \ref{lem-modified-connection} in the case that \(\psi\) is an eigenspinor.
We can estimate the energy momentum tensor by \(|T|^2\geq 2\lambda^2|\psi|^4\).
Thus, from \eqref{eigenvalue-modified-connection} we obtain
\begin{equation}
\lambda^2\geq K-2\Delta\log|\psi|,
\end{equation}
where \(K\) denotes the Gaussian curvature of \(M\). Note that the last term on the right hand side of \eqref{eigenvalue-modified-connection}
vanishes for \(\psi\) being an eigenspinor.
Integrating over the surface \(M\) and using \eqref{laplace-log} completes the proof of Theorem \ref{main-result}.

\section{Nodal sets of solutions of similar equations}
In this section we are concerned with the nodal set of \emph{harmonic spinors}, 
\emph{twistor spinors} and solutions of a non-linear Dirac equation.
Most of the results presented in this section are well-known in the literature.

\begin{Lem}[Bochner formula]
For an arbitrary spinor \(\psi\in\Gamma(\Sigma M)\) the following Bochner formula holds
\begin{equation}
\label{bochner}
\Delta\log|\psi|=\frac{K}{2}-\frac{\langle\psi,D^2\psi\rangle}{|\psi|^2}+\frac{|\nabla\psi|^2}{|\psi|^2}-\frac{1}{2}\frac{|d|\psi|^2|^2}{|\psi|^4},
\end{equation}
where \(K\) is the Gaussian curvature of the surface.
\end{Lem}
\begin{proof}
This follows by a direct calculation using \eqref{schroedinger-lichnerowicz}.
\end{proof}

\begin{Prop}
\label{prop-zeroset-harmonic}
Let \((M,g)\) be a closed Riemannian spin surface. Moreover, assume that \(\psi\in\Gamma(\Sigma^+M)\)
is a harmonic spinor, that is \(D\psi=0\). Then the following formula holds
\begin{equation}
N_0=-\frac{\chi(M)}{2},
\end{equation}
where \(N_0\) is the sum of the order of the zeros of \(\psi\) and \(\chi(M)\) denotes
the Euler characteristic of the surface.
\end{Prop}

\begin{proof}
Since \(\Sigma^+M\) is a complex line bundle, we can write \(\nabla_{e_j}\psi=f_j\psi\)
for some complex-valued function \(f_j\) away from its zero-set.
By a direct calculation we obtain the following identities
\begin{align}
\label{identity-a}\big|d|\psi|^2\big|^2=&4\sum_{j=1}^2|\operatorname{Re}f_j|^2|\psi|^4, \\
\label{identity-b}|\nabla\psi|^2=&\sum_{j=1}^2|f_j|^2|\psi|^2, \\
\label{identity-c}|D\psi|^2=&\big((\operatorname{Re}f_1+\operatorname{Im}f_2)^2+(\operatorname{Re}f_2+\operatorname{Im}f_1)^2\big)|\psi|^2.
\end{align}
By assumption \(D\psi=0\), consequently we get from \eqref{identity-c} that \(\operatorname{Re}f_1=-\operatorname{Im}f_2\)
and \(\operatorname{Re}f_2=-\operatorname{Im}f_1\). Inserting this into both \eqref{identity-a} and \eqref{identity-b}
yields
\begin{align}
\label{harmonic-spinor-derivative-identity}
\big|d|\psi|^2\big|^2=2|\psi|^2|\nabla\psi|^2.
\end{align}
Making use of the Bochner formula \eqref{bochner} and the assumption \(D\psi=0\) we obtain
\[
\Delta\log|\psi|=\frac{K}{2}.
\]
Integrating over the surface \(M\) and applying \eqref{laplace-log} then proves the assertion.
\end{proof}

\begin{Bem}
The same formula holds true if \(\psi\in\Gamma(\Sigma^-M)\) solves \(D\psi=0\).
Note that from a holomorphic perspective the statement of Proposition \ref{prop-zeroset-harmonic} is 
a standard fact about the zero divisor of a holomorphic section of a holomorphic line bundle, see \cite[p.144, Example 1]{MR1288523}.
\end{Bem}

The last Proposition is a special case of a structure Theorem for Dirac-harmonic maps between surfaces, which was proven in \cite{MR2496649}, Theorem 4.2. 
In addition, the statement was also proven using the Poincaré-Hopf Index Theorem in \cite{MR3210753}, Theorem 4.12.

Now, we turn to the analysis of \emph{twistor spinors}, for more details on them see the book \cite{MR1164864}.

\begin{Dfn}
Let \((M,g)\) be a \(n\)-dimensional Riemannian spin manifold. A twistor spinor is a section \(\psi\) of \(\Sigma M\)
satisfying
\begin{equation}
\label{twistor-spinor}
P\psi=0,
\end{equation}
where \(P_X\psi:=\nabla_X\psi+\frac{1}{n}X\cdot D\psi \) for every \(X\in TM\).
\end{Dfn}

For twistor spinors we will prove the following 
\begin{Prop}
Let \((M,g)\) be a closed Riemannian spin surface and let \(\psi\in\Gamma(\Sigma^+ M)\) be a twistor spinor.
Then we have
\begin{equation}
N(\psi)=\frac{1}{2}\chi(M),
\end{equation}
where \(N(\psi)\) is the sum of the order of the zeros of \(\psi\).
The same formula also holds for \(\psi\in\Gamma(\Sigma^- M)\).
\end{Prop}
\begin{proof}
We follow a similar approach as in the proof of the last Proposition.
First, if \(\psi\) is a twistor spinor on a \(n-\)dimensional Riemannian manifold, then it satisfies
\[
D^2\psi=R\frac{n}{4(n-1)}\psi,
\]
see for example \cite{MR1164864}, p.24.
Moreover, we can again write \(\nabla_{e_j}\psi=f_j\psi\) for some complex-valued function \(f_j\) away from its zero-set.
Since \(\psi\) is a twistor spinor, it satisfies
\begin{align}
\label{twistor-triangle}
|\nabla\psi|^2=\frac{1}{2}|D\psi|^2.
\end{align}
Combining \eqref{twistor-triangle} with \eqref{identity-b} and \eqref{identity-c} we obtain
\[
\big(|f_1|^2+|f_2|^2-\frac{1}{2}(\operatorname{Re}f_1+\operatorname{Im}f_2)^2-\frac{1}{2}(\operatorname{Re}f_2+\operatorname{Im}f_1)^2\big)|\psi|^2=0,
\]
which leads to
\[
\big((\operatorname{Re}f_1-\operatorname{Im}f_2)^2+(\operatorname{Re}f_2-\operatorname{Im}f_1)^2\big)|\psi|^2=0.
\]
Thus, we get \(\operatorname{Re}f_1=\operatorname{Im}f_2\) and \(\operatorname{Re}f_2=\operatorname{Im}f_1\)
such that \eqref{harmonic-spinor-derivative-identity} still holds.
In the end, via the Bochner formula (\ref{bochner}), we find
\begin{equation*}
\Delta\log|\psi|=-\frac{1}{2}K.
\end{equation*}
Integrating this equation using \eqref{laplace-log} and the Gauß-Bonnet theorem completes the proof.
\end{proof}

\begin{Cor}
The only closed surface admitting non-trivial twistor spinors with zeros is \(S^2\).
Moreover, there do not exist non-trivial twistor spinors on surfaces with negative Euler characteristic.
\end{Cor}

The spinorial Weierstraß representation for CMC-surfaces in \(\R^3\) is governed by the equation
\begin{equation}
\label{weierstrass}
D\psi=\mu|\psi|^2\psi,\qquad \|\psi\|_{L^4}=1,
\end{equation}
where \(\mu\) denotes some real constant. The condition on the \(L^4\) norm of \(\psi\) 
normalizes the volume of the surface to be equal to 1. We refer to \cite{MR2550205} and references therein for more details.
We can give an estimate on the nodal set of solutions of \eqref{weierstrass},
which was already proven in \cite{MR2550205}.
\begin{Prop}
The zero set of solutions of \eqref{weierstrass} can be estimated as 
\begin{equation}
N(\psi)\leq\frac{\mu^2}{4\pi}-\frac{\chi(M)}{2}.
\end{equation}
\end{Prop}
\begin{proof}
Again, by the main result of \cite{MR1714341} we know that the nodal set of solutions of \eqref{weierstrass} is discrete.
Moreover, the last term on the right hand side of \eqref{eigenvalue-modified-connection} vanishes when \(\psi\) 
solves \eqref{weierstrass}.
We compute
\begin{align*}
\langle\psi,D^2\psi\rangle=&\langle\psi,D(\mu|\psi|^2\psi)\rangle 
=\mu\underbrace{\langle\psi,(\nabla|\psi|^2)\cdot\psi\rangle}_{=0}+\mu\langle\psi,|\psi|^2D\psi\rangle 
=\mu^2|\psi|^6,
\end{align*}
where we used the skew-symmetry of the Clifford multiplication and \eqref{weierstrass} in the last step.
Using that \(|T|^2\geq 2\mu^2|\psi|^8\) in this case we obtain from \eqref{eigenvalue-equation-improved}
\[
\mu^2|\psi|^4\geq K-2\Delta\log|\psi|
\]
and integrating over the surface \(M\) using \eqref{laplace-log} yields the assertion.
\end{proof}

\section{Applications}
In this section we give some applications of Theorem \ref{main-result}.
\begin{Bsp}
We can use \eqref{estimate-nodal-set} to estimate the nodal set of the eigenspinor belonging to the first eigenvalue
on \(S^2\) with the round metric.
The Dirac spectrum on \(M=S^n\) for \(n\geq 2\) equipped with the round metric is well-known (see for example \cite{MR1361548}):
The eigenvalues of the Dirac operator are
\[
\lambda_k=\pm(\frac{n}{2}+k),\qquad k\in\N.
\]
Thus, from \eqref{estimate-nodal-set} we find that the zero set of an eigenspinor belonging to the first eigenvalue on \(S^2\) is
empty. This result is also well-known: On the sphere with the round metric the eigenspinors belonging to the eigenvalues \(\pm\frac{n}{2}\)
are Killing spinors with \(\alpha=\pm\frac{n}{2}\), see \cite{MR2509837}, Example A.1.3 and references therein. Thus, their nodal set
is empty. This statement can also be obtained by an explicit calculation:
In \cite{hermann-phd}, Chapter 4, it is shown that on the sphere \(S^n\) with the round metric
the eigenspinors corresponding to the eigenvalues \(\lambda=\pm\frac{n}{2}\)
are nowhere zero. 
\end{Bsp}

We can generalize this statement to the case of an arbitrary metric with positive curvature.
Recall the following upper bound for the first eigenvalue from \cite{MR1095763}, Corollary 1:
\begin{Lem}
Let \((M^{2m},g)\) be a compact even-dimensional spin manifold of positive
sectional curvature \(0\leq K\) and let \(K_{max}\) be the maximum of \(K\).
Then the first eigenvalue of the Dirac operator \(D\) is bounded by
\begin{equation}
\label{baum-inequal}
\lambda_1^2\leq 2^{n-2}\frac{n}{2}\max_M K^M.
\end{equation}
\end{Lem}

\begin{Lem}
Assume that \(M=S^2\) with a metric of positive curvature. 
Then we can estimate the nodal set of the eigenspinor belonging to the first eigenvalue on \(S^2\) as
\begin{equation}
N_1\leq \frac{\vol(S^2,g)}{4\pi}\max_{S^2}K-1.
\end{equation}
\end{Lem}

\begin{proof}
Applying \eqref{baum-inequal} on \(S^2\) 
\[
\lambda_1^2\leq\max_{S^2}K
\]
and using \eqref{estimate-nodal-set} proves the result.
\end{proof}

We can also give an upper bound for the nodal set on hyperbolic surfaces.
In this case we will make use of an inequality given by Lott \cite{MR860754}:
\begin{Satz}
Let \((M^n,g)\) be an \(n\)-dimensional closed spin manifold with \(n\geq 2\).
Then for any conformal class \([g]\) on \(M^n\), there exists \(b[g]>0\)
such that
\begin{equation}
\label{lott-inequality}
\lambda_1^2\leq b[g]\sup_M(-K)
\end{equation}
for any \(g\in[g]\) with \(R_g<0\).
\end{Satz}

\begin{Prop}
Let \(M\) be a closed surface with \(\chi(M)<0\). Then the nodal set of the eigenspinor corresponding
to the first eigenvalue of \(D\) can be estimated as
\begin{equation}
N_1\leq \frac{\vol(M,g)b[g]}{4\pi}\sup_M(-K)-\frac{\chi(M)}{2}. 
\end{equation}
\end{Prop}
\begin{proof}
This follows directly from \eqref{eigenvalue-equation-improved} and \eqref{lott-inequality}.
\end{proof}

The Willmore energy of a surface \(M\subset \R^3\) is defined by
\begin{equation*}
W(M):=\int_MH^2dM,
\end{equation*}
where \(H\) denotes the mean curvature of \(M\) in \(\R^3\).
It was shown in \cite{MR1651379} that the Willmore energy \(W(M)\) can be estimated 
with the help of the spectrum of the Dirac operator via
\[
W(M)\geq\lambda^2\vol(M,g).
\]
Using \eqref{eigenvalue-equation-improved} this directly implies

\begin{Prop}
We can estimate the Willmore energy as
\begin{equation}
W(M)\geq 2\pi\chi(M)+4\pi N(\psi).
\end{equation}
\end{Prop}

Using the Weyl-asymptotic for linear elliptic operators we can give an estimate on the nodal set
of eigenspinors for large eigenvalues. To this end we order the eigenvalues of \(D\)
by increasing absolute values, that is \(|\lambda_1|\leq|\lambda_2|\leq\ldots\).

\begin{Prop}
Let \(M\) be a closed surface and let \(\psi\in\Gamma(\Sigma M)\) be an eigenspinor of \(D\) with eigenvalue \(\lambda_k\).
For large values of \(k\) we find the following estimate
\begin{equation}
N(\psi)\leq \frac{4\pi^3}{\vol(M,g)}k^2-\frac{\chi(M)}{2}.
\end{equation}
\end{Prop}
\begin{proof}
Recall the Weyl asymptotic for large eigenvalues of a linear elliptic differential operator
\[
(\lambda_k)^\frac{n}{2}\sim\frac{(2\pi)^nk}{\vol(M,g)\omega_n},\qquad \omega_n=\frac{n^{n/2}}{\Gamma(\frac{n}{2}+1)}.
\]
The result follows by combining the Weyl asymptotic with \eqref{estimate-nodal-set}.
\end{proof}

\par\medskip
\textbf{Acknowledgements:}
The author would like to thank Andreas Hermann for several useful comments on a first draft of this note
and the referee whose remarks helped to significantly improve the content of the article.
In addition, the author gratefully acknowledges the support of the Austrian Science Fund (FWF)
through the START-Project Y963-N35 of Michael Eichmair.
\bibliographystyle{plain}
\bibliography{mybib}

\begin{thebibliography}{10}

\bibitem{MR2550205}
Bernd Ammann.
\newblock The smallest {D}irac eigenvalue in a spin-conformal class and cmc
  immersions.
\newblock {\em Comm. Anal. Geom.}, 17(3):429--479, 2009.

\bibitem{MR1900319}
Bernd Ammann and Christian B{\"a}r.
\newblock Dirac eigenvalue estimates on surfaces.
\newblock {\em Math. Z.}, 240(2):423--449, 2002.

\bibitem{MR1162671}
Christian B{\"a}r.
\newblock Lower eigenvalue estimates for {D}irac operators.
\newblock {\em Math. Ann.}, 293(1):39--46, 1992.

\bibitem{MR1361548}
Christian B{\"a}r.
\newblock The {D}irac operator on space forms of positive curvature.
\newblock {\em J. Math. Soc. Japan}, 48(1):69--83, 1996.

\bibitem{MR1651379}
Christian B{\"a}r.
\newblock Extrinsic bounds for eigenvalues of the {D}irac operator.
\newblock {\em Ann. Global Anal. Geom.}, 16(6):573--596, 1998.

\bibitem{MR1714341}
Christian B{\"a}r.
\newblock Zero sets of solutions to semilinear elliptic systems of first order.
\newblock {\em Invent. Math.}, 138(1):183--202, 1999.

\bibitem{MR1186015}
Christian B{\"a}r and Paul Schmutz.
\newblock Harmonic spinors on {R}iemann surfaces.
\newblock {\em Ann. Global Anal. Geom.}, 10(3):263--273, 1992.

\bibitem{MR1095763}
Helga Baum.
\newblock An upper bound for the first eigenvalue of the {D}irac operator on
  compact spin manifolds.
\newblock {\em Math. Z.}, 206(3):409--422, 1991.

\bibitem{MR1164864}
Helga Baum, Thomas Friedrich, Ralf Grunewald, and Ines Kath.
\newblock {\em Twistors and {K}illing spinors on {R}iemannian manifolds},
  volume 124 of {\em Teubner-Texte zur Mathematik [Teubner Texts in
  Mathematics]}.
\newblock B. G. Teubner Verlagsgesellschaft mbH, Stuttgart, 1991.
\newblock With German, French and Russian summaries.

\bibitem{MR600828}
Th. Friedrich.
\newblock Der erste {E}igenwert des {D}irac-{O}perators einer kompakten,
  {R}iemannschen {M}annigfaltigkeit nichtnegativer {S}kalarkr\"ummung.
\newblock {\em Math. Nachr.}, 97:117--146, 1980.

\bibitem{MR750754}
Th. Friedrich.
\newblock Zur {A}bh\"angigkeit des {D}irac-operators von der {S}pin-{S}truktur.
\newblock {\em Colloq. Math.}, 48(1):57--62, 1984.

\bibitem{MR1476425}
Thomas Friedrich.
\newblock {\em Dirac-{O}peratoren in der {R}iemannschen {G}eometrie}.
\newblock Advanced Lectures in Mathematics. Friedr. Vieweg \& Sohn,
  Braunschweig, 1997.
\newblock Mit einem Ausblick auf die Seiberg-Witten-Theorie. [With an outlook
  on Seiberg-Witten theory].

\bibitem{MR1806438}
Thomas Friedrich and Eui~Chul Kim.
\newblock Some remarks on the {H}ijazi inequality and generalizations of the
  {K}illing equation for spinors.
\newblock {\em J. Geom. Phys.}, 37(1-2):1--14, 2001.

\bibitem{MR2509837}
Nicolas Ginoux.
\newblock {\em The {D}irac spectrum}, volume 1976 of {\em Lecture Notes in
  Mathematics}.
\newblock Springer-Verlag, Berlin, 2009.

\bibitem{MR1288523}
Phillip Griffiths and Joseph Harris.
\newblock {\em Principles of algebraic geometry}.
\newblock Wiley Classics Library. John Wiley $\&$ Sons, Inc., New York, 1994.

\bibitem{hermann-phd}
Andreas Hermann.
\newblock Dirac eigenspinors for generic metrics.
\newblock {\em PhD-thesis}, 2012.

\bibitem{MR3210753}
Andreas Hermann.
\newblock Zero sets of eigenspinors for generic metrics.
\newblock {\em Comm. Anal. Geom.}, 22(2):177--218, 2014.

\bibitem{MR0358873}
Nigel Hitchin.
\newblock Harmonic spinors.
\newblock {\em Advances in Math.}, 14:1--55, 1974.

\bibitem{MR1738150}
Eui~Chul Kim and Thomas Friedrich.
\newblock The {E}instein-{D}irac equation on {R}iemannian spin manifolds.
\newblock {\em J. Geom. Phys.}, 33(1-2):128--172, 2000.

\bibitem{MR1031992}
H.~Blaine Lawson, Jr. and Marie-Louise Michelsohn.
\newblock {\em Spin geometry}, volume~38 of {\em Princeton Mathematical
  Series}.
\newblock Princeton University Press, Princeton, NJ, 1989.

\bibitem{MR860754}
John Lott.
\newblock Eigenvalue bounds for the {D}irac operator.
\newblock {\em Pacific J. Math.}, 125(1):117--126, 1986.

\bibitem{MR1474501}
R.~Schoen and S.~T. Yau.
\newblock {\em Lectures on harmonic maps}.
\newblock Conference Proceedings and Lecture Notes in Geometry and Topology,
  II. International Press, Cambridge, MA, 1997.

\bibitem{MR2496649}
Ling Yang.
\newblock A structure theorem of {D}irac-harmonic maps between spheres.
\newblock {\em Calc. Var. Partial Differential Equations}, 35(4):409--420,
  2009.

\end{thebibliography}
\end{document}